\documentclass[preprint,12pt]{elsarticle}

\usepackage{amssymb}
\usepackage[cp1251]{inputenc}
\usepackage{amsfonts, amssymb, amsthm, amsmath}
\usepackage{mathtext}
\usepackage{longtable}

\usepackage[english]{babel}

\renewcommand{\r}{\mathbb R}

\renewcommand{\le}{\leqslant}
\renewcommand{\ge}{\geqslant}

\renewcommand{\phi}{\varphi}

\newcommand{\eqd}{\stackrel{d}{=}}

\newcommand{\vX}{\emph{\textbf{X}}}
\newcommand{\vY}{\emph{\textbf{Y}}}
\newcommand{\vZ}{\emph{\textbf{Z}}}
\newcommand{\vS}{\emph{\textbf{S}}}
\newcommand{\vt}{\emph{\textbf{t}}}
\newcommand{\va}{\emph{\textbf{a}}}
\newcommand{\vc}{\emph{\textbf{c}}}
\newcommand{\vV}{\emph{\textbf{V}}}
\newcommand{\vW}{\emph{\textbf{W}}}
\newcommand{\vx}{\emph{\textbf{x}}}
\newcommand{\vy}{\emph{\textbf{y}}}

\newcommand{\vT}{\emph{\textbf{T}}}
\newcommand{\vw}{\emph{\textbf{w}}}
\newcommand{\vv}{\emph{\textbf{v}}}
\newcommand{\vI}{\emph{\textbf{I}}}

\newcommand{\bO}{\mathbf{0}}
\newcommand{\vb}{\emph{\textbf{b}}}

\journal{Journal of Statistical Planning and Inference}

\begin{document}

\begin{frontmatter}



\title{{On convergence of the distributions of statistics with
random sample sizes to 
normal variance-mean mixtures}}
\author[vyk]{V. Yu. Korolev}
\author[zai]{A. I. Zeifman}
\address[vyk]{Faculty of Computational Mathematics and
Cybernetics, Lomonosov Moscow State University; Institute of
Informatics Problems, Russian Academy of Sciences; e-mail victoryukorolev@yandex.ru}
\address[zai]{Vologda State University;
Institute of Informatics Problems, RAS; Institute of Socio-Economic
Development of Territories, RAS; e-mail a$\_$zeifman@mail.ru}

\begin{abstract}
We prove a general transfer theorem for multivariate random
sequences with independent random indexes in the double array limit
setting. We also prove its partial inverse providing necessary and
sufficient conditions for the convergence of randomly indexed random
sequences. Special attention is paid to the case where the elements
of the basic double array are formed as statistics constructed from
samples with random sizes. Under rather natural conditions we prove
the theorem on convergence of the distributions of such statistics
to multivariate normal variance-mean mixtures and, in particular, to
multivariate generalized hyperbolic laws.
\end{abstract}

\begin{keyword}
random sequence \sep random index \sep transfer theorem \sep
samples with random sizes \sep normal variance-mean mixture
\MSC[2010] 60F05 \sep 60G50 \sep 62E20
\end{keyword}

\end{frontmatter}

\section{Introduction}

In classical problems of mathematical statistics, the size of the
available sample, i. e., the number of available observations, is
traditionally assumed to be deterministic. In the asymptotic
settings it plays the role of infinitely increasing {\it known}
parameter. At the same time, in practice very often the data to be
analyzed is collected or registered during a certain period of time
and the flow of informative events each of which brings a next
observation forms a random point process, so that the number of
available observations is unknown till the end of the process of
their registration and also must be treated as a (random)
observation. For example, this is so in insurance statistics where
during different accounting periods different numbers of insurance
events (insurance claims and/or insurance contracts) occur; in
medical statistics where the number of patients with a certain
disease varies from month to month due to seasonal factors or from
year to year due to some epidemic reasons; in quality control where
the number of failed items differs from lot to lot; in
high-frequency financial statistics where the number of events in a
limit order book during a time unit essentially depends on the
intensity of order flows, etc. In these cases the number of
available observations as well as the observations themselves are
unknown beforehand and should be treated as random. Therefore it is
quite reasonable to study the asymptotic behavior of general
statistics constructed from samples with random sizes for the
purpose of construction of suitable and reasonable asymptotic
approximations. As this is so, to obtain non-trivial asymptotic
distributions in limit theorems of probability theory and
mathematical statistics, an appropriate centering and normalization
of random variables and vectors under consideration must be used. It
should be especially noted that to obtain reasonable approximation
to the distribution of the basic statistics, both centering and
normalizing values should be non-random. Otherwise the approximate
distribution becomes random itself and, for example, the problem of
evaluation of quantiles or significance levels becomes senseless.

In asymptotic settings, statistics constructed from samples with
random sizes are special cases of random sequences with random
indices. The randomness of indices usually leads to that the limit
distributions for the corresponding random sequences are
heavy-tailed even in the situations where the distributions of
non-randomly indexed random sequences are asymptotically normal see,
e. g., \cite{BeningKorolev2002, BeningKorolev2003,
GnedenkoKorolev1996}. For example, if a statistic which is
asymptotically normal in the traditional sense, is constructed on
the basis of a sample with random size having negative binomial
distribution, then instead of the expected normal law, the Student
distribution appears as an asymptotic law for this statistic.

The literature on random sequences with random indexes is extensive,
see, e. g., the references above and the references therein.

Although the mathematical theory of random sequences with random
indexes is well-developed, there still remain some unsolved
problems. For example, convenient conditions for the convergence of
the distributions of general statistics constructed from samples
with random sizes to normal variance-mean mixtures have not been
found yet. The desire to fill this theoretical gap is quite natural.
Another motivation for this research is practical and is as follows.
In applied probability there is a convention, apparently going back
to the book \cite{GnedenkoKolmogorov1954}, according to which a
model distribution is reasonable and/or justified enough only if it
is an {\it asymptotic approximation}, that is, there exist a more or
less simple limit setting and the corresponding limit theorem in
which the model under consideration is a limit distribution. The
existence of such a setting may bring a deeper insight into the
phenomena under consideration than just fitting a more or less
convenient model.

General normal variance-mean mixtures are examples of such
convenient models widely used to describe observed statistical
regularities in many fields. In particular, in 1977--78
O.\,Barndorff-Nielsen \cite{BN1977}, \cite{BN1978} introduced the
class of {\it generalized hyperbolic distributions} as a class of
special univariate variance-mean mixtures of normal laws in which
the mixing is carried out in one parameter since location and scale
parameters of the mixed normal distribution are directly linked. The
range of applications of generalized hyperbolic distributions varies
from the theory of turbulence or particle size description to
financial mathematics, see \cite{BN1982}. Multivariate generalized
hyperbolic distributions were introduced in the seminal paper
\cite{BN1977} mentioned above as a natural generalization of the
univariate case. They were further investigated in
\cite{Blaesild1981} and \cite{BlaesildJensen1981}. It is a
convention to explain such a good adequacy of generalized hyperbolic
models by that they possess many parameters to be suitably adjusted.
But actually, it would be considerably more reasonable to explain
this phenomenon by limit theorems yielding the possibility of the
use of generalized hyperbolic distributions as convenient {\it
asymptotic} approximations.

The main results presented in this paper deal with the description
of conditions which provide the convergence of the distributions of
statistics constructed from samples with random sizes to
multivariate normal variance-mean mixtures, in particular, to
multivariate generalized hyperbolic laws.
The conditions presented below are formulated in terms of the
asymptotic behavior of random sample sizes and have the `if and only
if' form. This circumstance proved to be very promising and
constructive. For example, in
\cite{KorolevChertokKorchaginZeifman2015} a problem of construction
of suitable approximations to the distribution of the so-called
order flow imbalance process in high-frequency trading systems was
considered. It was empirically shown that generalized hyperbolic
distributions are very likely models for that. But these
distributions are variance-mean mixtures with {\it one} mixing
parameter. By means of one-dimensional limit theorems for random
sums, this fact directly lead to theoretical understanding that the
intensities of buy and sell orders actually must be proportional to
{\it one and the same} random process reflecting general market
agitation. This theoretical inference concerning the flows
intensities then found its statistical proof, see
\cite{KorolevChertokKorchaginZeifman2015}. In other words, the `if
and only if' character of the presented conditions makes testing
goodness-of-fit of financial data with generalized hyperbolic models
{\it equivalent} to testing goodness-of-fit of the corresponding
flow intensities (i. e., volatilities) with generalized inverse
Gaussian models, which is much simpler.

In the present paper both structural and multivariate
generalizations of some results of
\cite{KorolevChertokKorchaginZeifman2015} to general statistics are
presented. The paper is organized as follows. Basic notation is
introduced in Section 2. Here an auxiliary result on the asymptotic
rapprochement of the distributions of randomly indexed random
sequences with special scale-location mixtures is proved. In Section
3 of the present paper we prove a general transfer theorem for
random sequences with independent random indexes in the double array
limit setting. We also prove its partial inverse providing necessary
and sufficient conditions for the convergence of randomly indexed
random indexes. Following the lines of \cite{Korolev1993}, we first
formulate a general result improving some results of
\cite{Korolev1993,BeningKorolev2002} by removing some superfluous
assumptions and relaxing some conditions. Special attention is paid
to the case where the elements of the basic double array are formed
as statistics constructed from samples with random sizes. This case
is considered in Section 4 where a theorem establishing necessary
and sufficient conditions of convergence to general multivariate
variance-mean mixtures is proved. As a corollary, in Section 5 we
deduce a criterion of convergence of the distributions of statistics
constructed from samples with random sizes to multivariate
generalized hyperbolic distributions.

\section{Notation. Auxiliary results}

Let us introduce basic notations to be used throughout this paper.
Let $m\in\mathbb{N}$. The vectors
$\vx=(x^{(1)},\ldots,x^{(m)})^{\top}$ are elements of
$\mathbb{R}^m$, the superscript $^{\top}$ stands for the transpose
of a vector or matrix. The scalar product in $\mathbb{R}^m$ will be
denoted $\langle\cdot,\,\cdot\rangle$:
$\langle\vx,\,\vy\rangle=\vx^{\top}\vy=x^{(1)}y^{(1)}+\ldots+x^{(m)}y^{(m)}$.
As usual, the Euclidean norm of $\vx$ is
$\|\vx\|=\langle\vx,\,\vx\rangle^{1/2}$. If $A$ is a real-valued
$(m\times m)$-square matrix, then $\det(A)$ denotes the determinant
of $A$. The $(m\times m)$-identity matrix is denoted $\vI$. To
properly distinguish between the real number zero and the zero
vector, we write $0\in\mathbb{R}$ and
$\bO=(0,\ldots,0)^{\top}\in\mathbb{R}^m$. The notation
$N_{\va,\Sigma}$ will be used for the $m$-dimensional normal
distribution with mean vector $\va$ and covariance matrix $\Sigma$.
The distribution function of the one-dimensional standard normal
distribution will be denoted $\Phi(x)$,
$$
\Phi(x)=\frac{1}{\sqrt{2\pi}}\int_{-\infty}^{x}e^{-y^2/2}dy,\ \ \
x\in\mathbb{R}.
$$
Assume that all the random variables and vectors considered in this
paper are defined on one and the same probability space
$(\Omega,\,\mathfrak{F},\,{\sf P})$. The symbols $\mathfrak{B}_m$
and $\mathfrak{B}_+$ will denote the Borel sigma-algebras of subsets
of $\mathbb{R}_m$ and $\mathbb{R}_+\equiv[0,\infty)$, respectively.
In what follows the symbols $\eqd$ and $\Longrightarrow$ will denote
coincidence of distributions and weak convergence (convergence in
distribution). We will write $\mathcal{L}(\vX)$ to denote the
distribution of a random vector $\vX$. A family
$\{\vX_j\}_{j\in\mathbb{N}}$ of $\mathbb{R}^m$-valued random vectors
is said to be {\it weakly relatively compact}, if each sequence of
its elements contains a weakly convergent subsequence. In the
finite-dimensional case the weak relative compactness of a family
$\{\vX_j\}_{j\in\mathbb{N}}$ is equivalent to its {\it tightness}
$\lim_{R\to\infty}\sup_{n\in\mathbb{N}}{\sf P}(\|\vX_n\|>R)=0$ (see,
e. g., \cite{Loeve1977}).

Let $\{\vS_{n,k}=(S_{n,k}^{(1)},\ldots,S_{n,k}^{(m)})^{\top}\}$,
$n,k\in\mathbb{N},$ be a double array of $\mathbb{R}^m$-valued
random vectors. For $n,k\in\mathbb{N}$ let
$\va_{n,k}=(a_{n,k}^{(1)},\ldots,a_{n,k}^{(m)})^{\top}\in\mathbb{R}^m$
be non-random vectors and $b_{n,k}\in\mathbb{R}$ be real numbers
such that $b_{n,k}>0$. The purpose of the vectors $\va_{n,k}$ and
numbers $b_{n,k}$ is to provide weak relative compactness of the
family of the random vectors $\big\{\vY_{n,k}\equiv
b_{n,k}^{-1}\big(\vS_{n,k}-\va_{n,k}\big)\big\}_{n,k\in\mathbb{N}}$
in the cases where it is required.

Consider a family $\{N_n\}_{n\in\mathbb{N}}$ of nonnegative integer
random variables such that for each $n,k\in\mathbb{N}$ the random
variables $N_n$ and random vectors $\vS_{n,k}$ are independent.
Especially note that we do not assume the row-wise independence of
$\{\vS_{n,k}\}_{k\ge1}$. Let
$\vc_n=(c_n^{(1)},\ldots,c_n^{(m)})^{\top}\in\mathbb{R}^m$ be
non-random vectors and $d_n$ be real numbers, $n\in\mathbb{N}$, such
that $d_n>0$. Our aim is to study the asymptotic behavior of the
random vectors $\vZ_n\equiv d_n^{-1}\big(\vS_{n,N_n}-\vc_n\big)$ as
$n\to\infty$ and find rather simple conditions under which the limit
laws for $\vZ_n$ have the form of normal variance-mean mixtures. In
order to do so we first formulate a somewhat more general result
following the lines of \cite{Korolev1993}, removing superfluous
assumptions, relaxing the conditions and generalizing the results of
that paper.

The characteristic functions of the random vectors $\vY_{n,k}$ and
$\vZ_n$ will be denoted $h_{n,k}(\vt)$ and $f_n(\vt)$, respectively,
$\vt\in\mathbb{R}^m$. Let $\vY$ be an $\mathbb{R}^m$-valued random
vector whose characteristic function will be denoted $h(\vt)$,
$\vt\in\mathbb{R}^m$. Introduce the random variables
$U_n=d_n^{-1}b_{n,N_n}$. Let
$\vV_n=(V_n^{(1)},\ldots,V_n^{(m)})^{\top}$ where
$V_n^{(k)}=d_n^{-1}(a_{n,N_n}^{(k)}-c_n^{(k)})$ is the $k$th
component of the random vector $d_n^{-1}(\va_{n,N_n}-\vc_n)$. In
what follows by $\vW_n$ we will denote the $(m+1)$-dimensional
compound random vector
$\vW_n=(U_n,\vV_n^{\top})^{\top}=(U_n,\,V_n^{(1)},\ldots,V_n^{(m)})^{\top}$.

Fot $\vt\in\mathbb{R}^m$ consider the function
$$
g_n(\vt)\equiv {\sf E}h(U_n\vt)e^{i\langle
\vt,\vV_n\rangle}={\textstyle\sum_{k=1}^{\infty}}{\sf
P}(N_n=k)e^{i\langle
\vt,d_n^{-1}(\va_{n,k}-\vc_n)\rangle}h\big(d_n^{-1}b_{n,k}\vt\big).\eqno(1)
$$
It can be easily seen that
$g_n(\vt)$ is the characteristic function of the random vector
$U_n\cdot \vY+\vV_n$ where the random vector $\vY$ is independent of
the random vector $\vW_n$.

In the double-array limit setting considered in this paper, to
obtain non-trivial limit laws for $\vZ_n$ we require the following
additional {\it coherency condition}: for any $T\in(0,\infty)$
$$
\lim_{n\to\infty}{\sf E}\sup_{\|\vt\|\le
T}\big|h_{n,N_n}(\vt)-h(\vt)\big|=0.\eqno(2)
$$

\smallskip

{\sc Remark 1}. It can be easily verified that, since the values
under the expectation sign in (2) are nonnegative and bounded (by
two), then coherency condition (2) is equivalent to that
$\sup_{\|\vt\|\le T}\big|h_{n,N_n}(\vt)-h(\vt)\big|\longrightarrow
0$ in probability as $n\to\infty$.

\smallskip

A particular form of the coherency condition depends on the
structure of the statistic $\vS_{n,k}$. For example, if $\vS_{n,k}$
is a sum of independent random variables and $h(\vt)$ is the
standard normal characteristic function, then, as it was
demonstrated in \cite{GKZ2014}, coherency condition (2) turns into
the easily verifiable {\it random Lindeberg condition}, whereas the
latter is not only sufficient, but necessary as well for the
convergence of the distributions of random sums of independent not
necessarily identically distributed random variables, see, e. g.,
\cite{KruglovZhangBo2001}.

\smallskip

{\sc Lemma 1.} {\it Let the family of random variables
$\{U_n\}_{n\in\mathbb{N}}$ be weakly relatively compact. Assume that
coherency condition $(2)$ holds. Then for any} $\vt\in\mathbb{R}^m$
{\it we have}
$$
\lim_{n\to\infty}|f_n(\vt)-g_n(\vt)|=0.\eqno(3)
$$

\smallskip

{\sc Proof}. Let $\gamma\in(0,\infty)$ be a real number to be
specified later. Denote $K_{1,n}\equiv K_{1,n}(\gamma)=\{k:\,
b_{n,k}\le\gamma d_n\}$, $K_{2,n}\equiv K_{2,n}(\gamma)=\{k:\,
b_{n,k}>\gamma d_n\}$. If $\vt=0$, then the assertion of the lemma
is trivial. Fix an arbitrary $\vt\neq 0$. By the formula of total
probability we have \setcounter{equation}{3}
\begin{eqnarray}
|f_n(\vt)-g_n(\vt)|=\nonumber\\
=\Big|{\textstyle\sum\nolimits_{k=1}^{\infty}}{\sf
P}(N_n=k)e^{i\langle\vt,\,d_n^{-1}(\va_{n,k}-\vc_n)\rangle}
\big[h_{n,k}\big(d_n^{-1}b_{n,k}\vt\big)-h\big(d_n^{-1}b_{n,k}\vt\big)\big]\Big|\le\nonumber\\
\le{\textstyle\sum\nolimits_{k\in K_{1,n}}}{\sf
P}(N_n=k)\big|h_{n,k}\big(d_n^{-1}b_{n,k}\vt\big)-h\big(d_n^{-1}b_{n,k}\vt\big)\big|+\nonumber\\
+{\textstyle\sum\nolimits_{k\in K_{2,n}}}{\sf
P}(N_n=k)\big|h_{n,k}\big(d_n^{-1}b_{n,k}\vt\big)-h\big(d_n^{-1}b_{n,k}\vt\big)\big|\equiv
I_1+I_2.
\end{eqnarray}
Choose an arbitrary $\epsilon>0$.

First consider $I_2$. We obviously have
$$
I_2\le 2{\textstyle\sum\nolimits_{k\in K_{2,n}(\gamma)}}{\sf
P}(N_n=k)=2{\sf P}\big(U_n>\gamma\big).\eqno(5)
$$
The weak relative compactness of the family
$\{U_n\}_{n\in\mathbb{N}}$ implies the existence of a
$\gamma_1=\gamma_1(\epsilon)$ such that $\sup_n{\sf
P}\big(U_n>\gamma_1\big)<\epsilon$. Therefore, setting
$\gamma=\gamma_1$ from (5) we obtain
$$
I_2<\epsilon.\eqno(6)
$$
Now consider $I_1$ with $\gamma$ chosen above. If $k\in
K_{1,n}(\gamma)$, then
$\big\|d_n^{-1}b_{n,k}\vt\big\|\le\gamma\|\vt\|$ and we have
\begin{eqnarray}I_1\le{\textstyle\sum\nolimits_{k\in
K_{1,n}(\gamma)}}{\sf P}(N_n=k)\sup_{\|\vx\|\le\gamma
\|\vt\|}|h_{n,k}(\vx)-h(\vx)|\le\nonumber\\
\le {\sf E}\sup_{\|\vx\|\le\gamma
\|\vt\|}|f_{n,N_n}(\vx)-h(\vx)|.\nonumber
\end{eqnarray}
Therefore, coherency condition (2) implies that there exists a
number $n_0=n_0(\epsilon,\gamma)$ such that for all $n\ge n_0$
$$
I_1<\epsilon.\eqno(7)
$$
Unifying (4), (6) and (7) we obtain that $|f_n(t)-g_n(t)|<2\epsilon$
for $n\ge n_0$. The arbitrariness of $\epsilon$ proves (3). The
lemma is proved.

\smallskip

Lemma 1 makes it possible to use the distribution defined by the
characteristic function $g_n(\vt)$ (see (1)) as an {\it accompanying
asymptotic} approximation to the distribution of the random vector
$\vZ_n$. In order to obtain a {\it limit} approximation, in the next
section we formulate and prove the transfer theorem.

\section{General transfer theorem and its inversion. The structure of limit laws}

{\sc Theorem 1}. {\it Assume that coherency condition $(2)$ holds.
If there exist a random variable $U$ and an $m$-dimensional random
vector} $\vV$ {\it such that the distributions of the
$(m+1)$-dimensional random vectors} $\vW_n$ {\it converge to that of
the random vector} $\vW=(U,\,\vV^{\top})^{\top}:$
$$
\vW_n\Longrightarrow \vW\ \ \ (n\to\infty),\eqno(8)
$$
{\it then}
$$
\vZ_n\Longrightarrow \vZ\eqd U\cdot\vY+\vV\ \ \
(n\to\infty).\eqno(9)
$$
{\it where the random vectors} $\vY$ {\it and}
$\vW=(U,\,\vV^{\top})^{\top}$ {\it are independent.}

\smallskip

{\sc Proof}. Treating $\vt\in\mathbb{R}^m$ as a fixed parameter,
represent the function $g_n(\vt)$ as $g_n(\vt)={\sf
E}h(U_n\vt)\exp\{i\langle\vt,\,\vV_n\rangle\}\equiv{\sf
E}\phi_{\vt}(\vW_n)$. Since for each $\vt\in\mathbb{R}^m$ the
function $\phi_{\vt}(\vw)\equiv
h(u\vt)\exp\{i\langle\vt,\vv\rangle\}$,
$\vw=(u,\vv^{\top})^{\top}\in \mathbb{R}^{m+1}$, is bounded and
continuous in $\vw$, then by the definition of the weak convergence
we have
$$
\lim_{n\to\infty}{\sf E}\phi_{\vt}(\vW_n)={\sf
E}\phi_{\vt}(\vW).\eqno(10)
$$
Using the Fubini theorem it can be easily verified that the function
on the right-hand side of (10) is the characteristic function of the
random variable $U\cdot\vY+\vV$ where the $m$-dimensional random
vector $\vY$ is independent of the $(m+1)$-dimensional random vector
$\vW$. Now the statement of the theorem follows from lemma 1 by the
triangle inequality. The theorem is proved.

\smallskip

It is easy to see that relation (9) is equivalent to that the limit
law for normalized randomly indexed random vectors $\vZ_n$ is a
scale-location mixture of the distributions which are limiting for
normalized non-randomly indexed random vectors $\vY_{n,k}$. Among
all scale-location mixtures, {\it variance-mean mixtures} attract a
special interest. To be more precise, we should speak of {\it normal
variance-mean mixtures} which are defined in the following way.

An $\mathbb{R}^m$-valued random vector $\vX$ is said to have a
multivariate normal mean-variance mixture distribution if
$\vX\eqd\va+U\vb+\sqrt{U}A\vY$, where $\va,\vb\in\mathbb{R}^m$, $A$
is a real $(m\times m)$-matrix such that the matrix $\Sigma\equiv
AA^{\top}$ is positive definite, $\vY$ is a random vector with the
standard normal distribution $N_{0,\vI}$ and $U$ is a real-valued,
non-negative random variable independent of $\vY$. Equivalently, a
probability measure $F$ on $(\mathbb{R}^m,\,\mathfrak{B}_m)$ is said
to be a multivariate normal mean-variance mixture if
$$
F(d\vx)=\int_{0}^{\infty}N_{\vb+z\va,\,z\Sigma}(d\vx)G(dz),
$$
where the mixing distribution $G$ is a probability measure on
$(\mathbb{R}_+,\,\mathfrak{B}_+)$. In this case we will sometimes
briefly write $F = N_{\vb+z\va,\,z\Sigma}\circ G$.

Let us see how these mixtures can appear in the double-array setting
under consideration. Assume that the centering vectors $\va_{n,k}$
and $\vc_n$ are in some sense proportional to the scaling constants
$b_{n,k}$ and $d_n$. Namely, assume that there exist vectors
$\va_n\in\mathbb{R}^m$ and $\vb_n\in\mathbb{R}^m$ such that for all
$n,k\in\mathbb{N}$ we have $\va_{n,k}=d_n^{-1}b_{n,k}^2\va_n$,
$\vc_n=d_n\vb_n$, and there exist finite limits
$\va=\lim_{n\to\infty}\va_n$, $\vb=\lim_{n\to\infty}\vb_n$. Then
under condition (8)
$\vW_n=\big(U_n,\,(U_n^2\va_n+\vb_n)^{\top}\big)\vspace{-1mm}^{\top}\vspace{1mm}\Longrightarrow
\big(U,\,(U^2\va+\vb)^{\top}\big)\vspace{-1mm}^{\top}\vspace{1mm}$
$(n\to\infty)$, so that if in theorem 2 $\vY$ has the
$m$-dimensional normal distribution $N_{\bO,\,\Sigma}$, then the
limit law for $\vZ_n$ takes the form of the normal variance-mean
mixture $N_{\vb+z\va,\,z\Sigma}\circ G$ with $G$ being the
distribution of $U^2$.

In order to prove a result that is a partial inversion of theorem 1,
for fixed random vectors $\vZ$ and $\vY$ with the characteristic
functions $f(\vt)$ and $h(\vt)$ introduce the set
$\mathcal{W}(\vZ|\vY)$ containing all $(m+1)$-dimensional random
vectors $\vW=(U,\,\vV^{\top})^{\top}$ with $U\in\mathbb{R}$ and
$\vV\in\mathbb{R}^m$ such that the characteristic function $f(\vt)$
can be represented as
$$
f(\vt)={\sf E}h(U\vt)e^{i\langle\vt,\vV\rangle},\ \ \
\vt\in\mathbb{R}^m,\eqno(12)
$$
and ${\sf P}(U\ge 0)=1.$ Whatever random vectors $\vZ$ and $\vY$
are, the set $\mathcal{W}(\vZ|\vY)$ is always nonempty since it
trivially contains the vector $(0,\vZ^{\top})^{\top}$. It is easy to
see that representation (12) is equivalent to that $\vZ\eqd
U\vY+\vV$.

The set $\mathcal{W}(\vZ|\vY)$ may contain more that one element.
For example, if $\vY$ is the random vector with standard normal
distribution $N_{\bO,\vI}$ and $\vZ\eqd \vT_1-\vT_2$ where $\vT_1$
and $\vT_2$ are independent random vectors with independent
components having the same standard exponential distribution, then
along with the vector
$\big(0,(\vT_1-\vT_2)^{\top}\big)\vspace{-1mm}^{\top}\vspace{1mm}$
the set $\mathcal{W}(\vZ|\vY)$ contains the vector
$\big(\sqrt{U},\bO^{\top}\big)\vspace{-1mm}^{\top}\vspace{1mm}$
where $U$ is a random variable with the standard exponential
distribution. In this case $\vZ$ has the spherically symmetric
Laplace distribution.

Let $\Lambda(\vX_1,\,\vX_2)$ be any probability metric which
metrizes weak convergence in the space of $(m+1)$-dimensional random
vectors. An example of such a metric is the L{\'e}vy--Prokhorov
metric (see, e. g., \cite{Billingsley1968} or \cite{Zolotarev1997}).

\smallskip

{\sc Theorem 2}. {\it Let the family of random variables
$\{U_n\}_{n\in\mathbb{N}}$ be weakly relatively compact. Assume that
coherency condition $(2)$ holds. Then a random vector} $\vZ$ {\it
such that}
$$
\vZ_n\Longrightarrow \vZ\ \ \ (n\to\infty)\eqno(13)
$$
{\it with some} $\vc_n\in\mathbb{R}^m$ {\it exists if and only if
there exists a weakly relatively compact sequence of random vectors}
$\vW_n^*\equiv(U_n^*,\,(\vV_n^*)^{\top})^{\top}\in\mathcal{W}(\vZ|\vY)$,
$n\in\mathbb{N}$, {\it such that}
$$
\lim_{n\to\infty}\Lambda(\vW_n^*,\,\vW_n)=0.\eqno(14)
$$

\smallskip

{\sc Proof}. {\it ``Only if'' part}. Prove that the sequence
$\{\vV_n\}_{n\in\mathbb{N}}$ is weakly relatively compact. The
indicator function of a set $A$ will be denoted $\mathbb{I}(A)$. By
the formula of total probability for an arbitrary $R>0$ we have
\begin{eqnarray}
{\sf P}(\|\vV_n\|>R)={\textstyle\sum\nolimits_{k=1}^{\infty}}{\sf
P}(N_n=k)\mathbb{I}(\|d_n^{-1}(\va_{n,k}-\vc_n)\|>R)=\nonumber\\
= {\textstyle\sum\nolimits_{k=1}^{\infty}}{\sf P}(N_n=k){\sf
P}\big(\|d_n^{-1}(\vS_{n,k}-\vc_n)-d_n^{-1}b_{n,k}\cdot
b_{n,k}^{-1}(\vS_{n,k}-\va_{n,k})\|>R\big)\le\nonumber\\
\le {\sf
P}(2\|\vZ_n\|>R)+{\textstyle\sum\nolimits_{k=1}^{\infty}}{\sf
P}(N_n=k){\sf
P}(2d_n^{-1}b_{n,k}\cdot\|\vY_{n,k}\|>R)\equiv\nonumber\\
I_{1,n}(R)+I_{2,n}(R).\nonumber
\end{eqnarray}
First consider $I_{2,n}(R)$. Using the set $K_{2,n}=K_{2,n}(\gamma)$
introduced in the preceding section, for an arbitrary $\gamma>0$ we
have
\setcounter{equation}{14}
\begin{eqnarray}
I_{2,n}(R)={\textstyle\sum\nolimits_{k\in K_{2,n}}}{\sf
P}(N_n=k){\sf P}(2\|\vY_{n,k}\|>Rb_{n,k}^{-1}d_n)+\nonumber\\
+{\textstyle\sum\nolimits_{k\notin K_{2,n}}}{\sf P}(N_n=k){\sf
P}(2\|\vY_{n,k}\|>Rb_{n,k}^{-1}d_n)\le\nonumber\\
\le {\textstyle\sum\nolimits_{k\in K_{2,n}}}{\sf P}(N_n=k){\sf
P}(2\gamma\|\vY_{n,k}\|>R)+{\sf P}(U_n>\gamma)\le\nonumber\\
\le {\sf P}(2\gamma\|\vY_{n,N_n}\|>R)+{\sf P}(U_n>\gamma).
\end{eqnarray}
Fix an arbitrary $\epsilon>0$. Choose
$\gamma=\gamma(\epsilon)$ so that
$$
{\sf P}\big(U_n>\gamma(\epsilon)\big)<\epsilon\eqno(16)
$$
for all $n\in\mathbb{N}$. This is possible due to the weak relative
compactness of the family $\{U_n\}_{n\in\mathbb{N}}$. Now choose
$R'=R'(\epsilon)$ so that
$$
{\sf
P}(2\gamma(\epsilon)\|\vY_{n,N_n}\|>R'(\epsilon))<\epsilon.\eqno(17)
$$
This is possible due to the weak relative compactness of the family
$\{\vY_{n,N_n}\}_{n\in\mathbb{N}}$ implied by coherency condition
(2). Thus, from (15), (16) and (17) we obtain
$$
I_{2,n}\big(R'(\epsilon)\big)<2\epsilon\eqno(18)
$$
for all $n\in\mathbb{N}$. Now consider $I_{1,n}(R)$. From (13) it
follows that there exists an $R''=R''(\epsilon)$ such that
$$
I_{1,n}\big(R''(\epsilon)\big)<\epsilon\eqno(19)
$$
for all $n\in\mathbb{N}$. From (18) and (19) it follows that if
$R>\max\big\{R',\,R''\big\}$, then $\sup_{n}{\sf
P}(\|\vV_n\|>R)<3\epsilon$ and by virtue of the arbitrariness of
$\epsilon>0$, the family $\{\vV_n\}_{n\in\mathbb{N}}$ is weakly
relatively compact. Hence, the family of random vectors
$\{\vW_n\}_{n\in\mathbb{N}}$ is weakly relatively compact.

Denote $\epsilon_n=\inf\{\Lambda(\vW_n,\,\vW):\
\vW\in\mathcal{W}(\vZ|\vY)\}$, $n=1,2,...$ Prove that $\epsilon_n\to
0$ as $n\to\infty$. Assume the contrary. In this case $\epsilon_n\ge
M$ for some $M>0$ and all $n$ from some subsequence $\mathcal{N}$ of
natural numbers. Choose a subsequence
$\mathcal{N}_1\subseteq\mathcal{N}$ so that the sequence of random
vectors $\{\vW_n\}_{n\in\mathcal{N}_1}$ weakly converges to some
random vector $\vW$. As this is so, for all $n\in\mathcal{N}_1$
large enough we will have $\Lambda(\vW_n,\,\vW)<M$. Applying theorem
1 to the sequence $\{\vW_n\}_{n\in\mathcal{N}_1}$ we make sure that
$\vW\in \mathcal{W}(\vZ|\vY)$ since condition (13) implies the
coincidence of the limits of all convergent subsequences of
$\{\vZ_n\}$. We arrive at the contradiction with the assumption that
$\epsilon_n>M$ for all $n\in\mathcal{N}_1$. Hence, $\epsilon_n\to 0$
as $n\to\infty$. For each $n\in\mathbb{N}$ choose a vector
$\vW_n^*\in\mathcal{W}(\vZ|\vY)$ such that
$\Lambda(\vW_n,\,\vW_n^*)\le\epsilon_n+{\textstyle\frac1n}$. The
sequence $\{\vW_n^*\}_{n\in\mathbb{N}}$ obviously satisfies
condition (14). Its weak relative compactness follows from (14) and
the weak relative compactness of the sequence
$\{\vW_n\}_{n\in\mathbb{N}}$ established above.

\smallskip

{\it ``If'' part}. Assume that the sequence
$\{\vZ_n\}_{n\in\mathbb{N}}$ does not converge weakly to $\vZ$ as
$n\to\infty$. In that case the inequality
$\Lambda\big((0,\,\vZ_n^{\top})^{\top},\,(0,\,\vZ^{\top})^{\top})\ge
M$ holds for some $M>0$ and all $n$ from some subsequence
$\mathcal{N}$ of natural numbers. Choose a subsequence
$\mathcal{N}_1\subseteq\mathcal{N}$ so that the sequence of random
vectors
$\{\vW_n^*=(U_n^*,(\vV_n^*)^{\top})^{\top}\}_{n\in\mathcal{N}_1}$
weakly converges to some random vector $\vW$. Repeating the
reasoning used to prove theorem 1 we make sure that $ {\sf
E}e^{i\langle\vt,\vZ\rangle}={\sf
E}h(U_n^*\vt)e^{i\langle\vt,\vV_n^*\rangle}\longrightarrow{\sf
E}h(U\vt)e^{i\langle\vt,\vV\rangle} $ as $n\to\infty$,
$n\in\mathcal{N}_1$, for any $t\in\mathbb{R}$ that is,
$\vW\in\mathcal{W}(\vZ|\vY)$. From the triangle inequality
$\Lambda(\vW_n,\,\vW)\le
\Lambda(\vW_n,\,\vW_n^*)+\Lambda(\vW_n^*,\,\vW)$ and condition (14)
it follows that $\Lambda(\vW_n,\,\vW)\to0$ as $n\to\infty$,
$n\in\mathcal{N}_1$. Apply theorem 1 to the double array
$\{\vY_{n,k}\}_{k\in\mathbb{N},\,n\in\mathcal{N}_1}$ and the
sequence $\{\vW_n\}_{n\in\mathcal{N}_1}$. As a result we obtain that
$\Lambda((0,\vZ_n^{\top})^{\top},(0,\vZ^{\top})^{\top})\to0$ as
$n\to\infty$, $n\in\mathcal{N}_1$, contradicting the assumption that
$\Lambda((0,\vZ_n^{\top})^{\top},(0,\vZ^{\top})^{\top})\ge M>0$ for
$n\in\mathcal{N}_1$. Thus, the theorem is completely proved.

\smallskip

{\sc Remark 2}. It should be noted that in \cite{Korolev1993} and
some subsequent papers a stronger and less convenient version of the
coherency condition was used. Furthermore, in \cite{Korolev1993} and
the subsequent papers the statements analogous to lemma 1 and
theorems 1 and 2 were proved under the additional assumption of the
weak relative compactness of the family
$\{\vY_{n,k}\}_{n,k\in\mathbb{N}}$.

\section{Limit theorems for statistics constructed from samples with random sizes}

Let $\{\vX_{n,j}\}_{j\ge1}$, $n\in\mathbb{N},$ be a double array of
row-wise independent not necessarily identically distributed random
vectors with values in $\mathbb{R}^r$, $r\in\mathbb{N}$. For
$n,k\in\mathbb{N}$ let
$\vT_{n,k}=\vT_{n,k}(\vX_{n,1},...,\vX_{n,k})$ be a statistic, i.e.,
a measurable function of $\vX_{n,1},...,\vX_{n,k}$ with values in
$\mathbb{R}^m$. For each $n\geq1$ we define a random vector
$\vT_{n,N_n}$ by setting $\vT_{n,N_n}(\omega)\equiv
\vT_{{n,N_n}(\omega)}(\vX_{n,1}(\omega),...,\vX_{n,N_n(\omega)}(\omega))$,
$\omega\in\Omega$.

Let $\theta_n$ be $\mathbb{R}^m$-valued vectors, $n\in\mathbb{N}$.
In this section we will assume that the random vectors $\vS_{n,k}$
have the form $\vS_{n,k}=\vT_{n,k}-\theta_n$, $n,k\in\mathbb{N}$.
Concerning the normalizing constants and vectors we will assume that
there exist $m$-dimensional vectors $\va$, $\va_n$, $\vb$, $\vb_n$
and positive numbers $\sigma_n$ such that
$$
\va_n\to\va,\ \ \ \vb_n\to\vb \ \ \ (n\to\infty)\eqno(20)
$$
and for all $n,k\in\mathbb{N}$
$$
b_{n,k}=(\sigma_n\sqrt{k})^{-1},\ \ d_n=(\sigma_n\sqrt{n})^{-1},\ \
\va_{n,k}=(\sigma_n k)^{-1}\sqrt{n}\va_n,\ \
\vc_n=(\sigma_n\sqrt{n})^{-1}\vb_n\eqno(21)
$$
so that
$$
\vY_{n,k}=\sigma_n\sqrt{k}(\vT_{n,k}-\theta_n)-\sqrt{n/k}\va_n\ \
\text{and}\ \ \vZ_n=\sigma_n\sqrt{n}(\vT_{n,N_n}-\theta_n)-\vb_n.
$$
As this is so, $\sigma_n^2\vI$ can be regarded as the asymptotic
variance of $\vT_{n,k}$ as $k\to\infty$ whereas the bias of
$\vT_{n,k}$ is $\sqrt{n}(k\sigma_n)^{-1}\va_n$.

Recall that the characteristic function of the normal distribution
in $\mathbb{R}^m$ with zero expectation and covariance matrix
$\Sigma$ is $\phi(\vt)=\exp\{-\frac12\vt^{\top}\Sigma\vt\}$,
$\vt\in\mathbb{R}^m$. In what follows we will assume that the
statistic $\vT_{n,k}$ is asymptotically normal in the following
sense: there exists a positive definite symmetric matrix $\Sigma$
such that for any $T\in(0,\infty)$
$$
\lim_{n\to\infty}{\sf E}\sup_{\|\vt\|\le
T}\big|h_{n,N_n}(\vt)-\exp\{-{\textstyle\frac12}\vt^{\top}\Sigma\vt\big|=0,\eqno(22)
$$
where $h_{n,k}(\vt)$ is the characteristic function of the random
vector $\vY_{n,k}$.

\smallskip

{\sc Theorem 3}. {\it Let the family of random variables
$\{n/N_n\}_{n\in\mathbb{N}}$ be weakly relatively compact, the
normalizing constants have the form $(21)$ and satisfy condition
$(20)$. Assume that the statistic} $\vT_{n,k}$ {\it is
asymptotically normal so that condition $(22)$ holds. Then a random
vector} $\vZ$ {\it such that}
$$
\sigma_n\sqrt{n}(\vT_{n,N_n}-\theta_n)-\vb_n\Longrightarrow Z\ \ \
(n\to\infty)
$$
{\it exists if and only if there exists a distribution function $G$
such that $G(0)=0$, the distribution $F$ of} $\vZ$ {\it has the
form} $F=N_{\vb+z\va,\,z\Sigma}\circ G$ {\it and}
$$
{\sf P}(n/N_n<x)\Longrightarrow G(x)\ \ \ (n\to\infty).\eqno(23)
$$

\smallskip

{\sc Proof}. We will deduce theorem 3 as a corollary of theorem 2.
First, notice that condition (22) is actually the coherency
condition (2) with $h(\vt)\equiv
\exp\{-\frac12\vt^{\top}\Sigma\vt\}$.

Second, notice that, obviously, each one-dimensional marginal
distribution of a multivariate normal variance-mean mixture is a
one-dimensional normal variance-mean mixture. Recently in
\cite{Korolev2013} it was proved that one-dimensional normal
variance-mean mixtures are identifiable, that is, if
$a\in\mathbb{R}$, $\sigma>0$, ${\sf P}(Y<x)\equiv\Phi(x)$ and $U_1$
and $U_2$ are two nonnegative random variables, then the identity
$$
{\sf E}\Phi\Big(\frac{x-aU_1}{\sigma\sqrt{U_1}}\Big)\equiv {\sf
E}\Phi\Big(\frac{x-aU_2}{\sigma\sqrt{U_2}}\Big)
$$
implies that $U_1\eqd U_2$. Hence it follows that the set
$\mathcal{W}(\vZ|\vY)$ contains at most one vector of the form
$\vW=(\sigma\sqrt{U},\,(U\va+\vb)^{\top})^{\top}$. This means that
in the case under consideration condition (14) reduces to (23). The
theorem is proved.

\smallskip

{\sc Remark 3}. Note that the statistics $\vT_{n,k}$ in the
coherency condition are centered, whereas if $\vb_n=\vb=\bO$ and
$\va\neq \bO$, then the limit distribution for statistics
constructed from samples with random sizes becomes skew unlike in
the classical situation, where the presence of the systematic bias
of the original statistic results in that the limit distribution
becomes just shifted. So, if the limit normal variance-mean mixture
is skew, then it can be suspected that the original statistics are
actually biased.

\smallskip

The class of normal variance-mean mixtures is very wide. For
example, it contains generalized hyperbolic laws with generalized
inverse Gaussian mixing distributions, in particular, $(a)$
symmetric and non-symmetric (skew) Student distributions (including
Cauchy distribution), to which there correspond inverse gamma mixing
distributions; $(b)$ variance gamma (VG) distributions) (including
symmetric and non-symmetric Laplace distributions), to which there
correspond gamma mixing distributions; $(c)$
normal$\backslash\!\backslash$inverse Gaussian (NIG) distributions
to which there correspond inverse Gaussian mixing distributions, and
many other types. Along with generalized hyperbolic laws, the class
of normal variance-mean mixtures contains symmetric strictly stable
laws with strictly stable mixing distributions concentrated on the
positive half-line, generalized exponential power distributions and
many other types. By variance-mean mixing many other initially
symmetric types represented as pure scale mixtures of normal laws
can be skewed, e. g., as it was done to obtain non-symmetric
exponential power distributions in \cite{GK2013}.

\section{Convergence to multivariate generalized hyperbolic
distributions}

Generalized hyperbolic distributions demonstrate exceptionally high
adequacy when they are used to describe statistical regularities in
the behavior of characteristics of various complex open systems, in
particular, turbulent systems and financial markets. There are
dozens of dozens of publications dealing with models based on
univariate and multivariate generalized hyperbolic distributions.
Just mention the canonic papers \cite{BN1978, BN1979, 
EberleinKeller1995, Prause1997, 
EberleinKellerPrause1998, BarndorffNielsen1998, EberleinPrause1998,
Eberlein1999, BNBlaesildSchmiegel2004}. Therefore below we will
concentrate our attention on limit theorems establishing the
convergence of the distributions of statistics constructed from
samples with random sizes to multivariate generalized hyperbolic
distributions.

In order to do so we should first recall the definition of the {\it
generalized inverse Gaussian distribution} $GIG_{\nu,\mu,\lambda}$
on $\mathfrak{B}_+$. The density of this distribution is denoted
$p_{GIG}(x;\nu,\mu,\lambda)$ and has the form
$$
p_{GIG}(x;\nu,\mu,\lambda)=\frac{\lambda^{\nu/2}}{2\mu^{\nu/2}K_{\nu}\big(\sqrt{\mu\lambda}\big)}\cdot
x^{\nu-1}\cdot\exp\Big\{-\frac12\Big(\frac{\mu}{x}+\lambda
x\Big)\Big\},\ \ \ x>0.
$$
Here $\nu\in\r$,
$$
\begin{array}{lll}
\mu>0, & \lambda\ge0, & \text{if }\nu<0,\vspace{1mm}\cr \mu>0, &
\lambda>0, & \text{if }\nu=0,\vspace{1mm}\cr \mu\ge0, & \lambda>0, &
\text{if }\nu>0,
\end{array}
$$
$K_{\nu}(z)$ is the modified Bessel function of the third kind with
index $\nu$,
$$
K_{\nu}(z)=\frac12\int_{0}^{\infty}y^{\nu-1}\exp\Big\{-\frac{z}{2}\Big(y+\frac1y\Big)\Big\}dy,\
\ \ \ z\in\mathbb{C},\ \mathrm{Re}\,z>0.
$$
According to \cite{Seshadri1997}, the generalized inverse Gaussian
distribution was introduced in 1946 by {\'E}tienne Halphen, who used
it to describe monthly volumes of water passing through
hydroelectric power stations. In the paper \cite{Seshadri1997}
generalized inverse Gaussian distribution was called the {\it
Halphen distribution}. In 1973 this distribution was re-discovered
by Herbert Sichel \cite{Sichel1973}, who used it as the mixing law
in special mixed Poisson distributions (the {\it Sichel
distributions}, see, e. g., \cite{KorolevBeningShorgin2011}) as
discrete distributions with heavy tails. In 1977 these distributions
were once more re-discovered by O. Barndorff-Nielsen \cite{BN1977,
BN1978}, who, in particular, used them to describe the particle size
distribution.

The class of generalized inverse Gaussian distributions is rather
rich and contains, in particular, both distributions with
exponentially decreasing tails (gamma-distribution ($\mu=0$,
$\nu>0$)), and distributions whose tails demonstrate power-type
behavior (inverse gamma-distribution ($\lambda=0$, $\nu<0$), inverse
Gaussian distribution ($\nu=-\frac12$) and its limit case as
$\lambda\to0$, the L{\'e}vy distribution (stable distribution with
the characteristic exponent equal to $\frac12$ and concentrated on
the nonnegative half-line, the distribution of the time for the
standard Wiener process to hit the unit level)).

In the final part of his seminal paper \cite{BN1977}, O.
Barndorff-Nielsen defined the class of multivariate {\it generalized
hyperbolic distributions} as the class of special normal
variance-mean mixtures. Namely, let $\Sigma$ be a positive definite
$(m\times m)$-matrix with $\det(\Sigma)$=1, $\va$ and $\vb$ be
$m$-dimensional vectors. Then the $m$-dimensional generalized
hyperbolic distribution $GH_{\nu,\mu,\lambda,\va,\vb,\Sigma}$ on
$\mathfrak{B}_m$ is defined as
$$
GH_{\nu,\mu,\alpha,\va,\vb,\Sigma}=\mathcal{N}_{\vb+z\Sigma\va,\,z\Sigma}\circ
GIG(\nu,\mu,\sqrt{\alpha^2-\langle\va,\Sigma\va\rangle}).
$$
Due to the restrictions imposed on the parameters of the generalized
inverse Gaussian distribution, the parameters of generalized
hyperbolic distribution must fit the conditions $\nu\in\mathbb{R}$,
$\alpha,\mu\in\mathbb{R}_+$ and
$$
\begin{array}{lll}
\mu>0, & 0\le\langle\va,\Sigma\va\rangle\le\alpha^2, & \text{if
}\nu<0,\vspace{1mm}\cr \mu>0, &
0\le\langle\va,\Sigma\va\rangle<\alpha^2, & \text{if
}\nu=0,\vspace{1mm}\cr \mu\ge0, &
0\le\langle\va,\Sigma\va\rangle<\alpha^2, & \text{if
}\nu>0,\vspace{1mm}
\end{array}
$$
The corresponding distribution density
$p_{GH}(\vx;\nu,\mu,\alpha,\va,\vb,\Sigma)$ has the form
$$
\begin{array}{lll}
p_{GH}(\vx;\nu,\mu,\alpha,\va,\vb,\Sigma)=\vspace{3mm}\\
{\displaystyle=\frac{(\alpha^2-\langle\va,\Sigma\va\rangle)^{\nu/2}}{(2\pi)^{m/2}\alpha^{\nu-m/2}
\mu^{\nu/2}K_{\nu}\big(\sqrt{\mu(\alpha^2-\langle\va,\Sigma\va\rangle)}\big)}\sqrt{(\langle\vx-\vb,\Sigma^{-1}(\vx-\vb)\rangle+\mu)^{\nu-m/2}}\times}
\vspace{3mm}\\
\times
K_{\nu-m/2}\big(\alpha\sqrt{\langle\vx-\vb,\Sigma^{-1}(\vx-\vb)\rangle+\mu}\big)\exp\{\langle\va,\vx-\vb\rangle\},\
\ \ \ \ \ \vx\in\mathbb{R}^m.
\end{array}
$$

\smallskip

{\sc Theorem 4.} {\it Let the family of random variables
$\{n/N_n\}_{n\in\mathbb{N}}$ be weakly relatively compact, the
normalizing constants have the form $(21)$ and satisfy condition
$(20)$ with some} $\va,\vb\in\mathbb{R}^m$. {\it Assume that the
statistic} $\vT_{n,k}$ {\it is asymptotically normal so that
condition $(22)$ holds with some symmetric positive definite matrix
$\Sigma$. Then the distribution of a statistic} $\vT_{n,N_n}$ {\it
constructed from the sample with random size $N_n$ weakly converges,
as $n\to\infty$, to an $m$-dimensional generalized hyperbolic
distribution}:
$$
\mathcal{L}\big(\sigma_n\sqrt{n}(\vT_{n,N_n}-\theta_n)-\vb_n\big)\Longrightarrow
GH_{\nu,\mu,\alpha,\Sigma^{-1}\va,\vb,\Sigma}
$$
{\it if and only if}
$$
\mathcal{L}\big(n^{-1}N_n\big)\Longrightarrow
GIG_{-\nu,\lambda,\mu}\eqno(24)
$$
{\it with} $\lambda=\sqrt{\alpha^2-\langle\va,\Sigma\va\rangle}$.

\smallskip

This theorem is a straightforward corollary of theorem 3 with the
account of a simply verifiable fact that if
$\mathcal{L}(\xi)=GIG_{\nu,\mu,\lambda}$, then
$\mathcal{L}(\xi^{-1})=GIG_{-\nu,\lambda,\mu}$.

\smallskip

Theorem 4 can serve as convenient explanation of the high adequacy
of generalized hyperbolic L{\'e}vy distributions as models of
statistical regularities in the behavior of stochastic systems.
Moreover, they directly link the mixing distribution in the
representation of a generalized hyperbolic distribution with the
random sample size which is determined by the intensity of the flow
of informative events generating the observations, see, e. g.,
\cite{KorolevChertokKorchaginZeifman2015}.

According to theorem 4, for example, to obtain the limit
multivariate asymmetric Student distribution for $\vT_{n,N_n}$ it is
necessary and sufficient that in (24) the mixing distribution is the
gamma distribution \cite{KorolevSokolov2012}. To obtain the
multivariate variance gamma limit distribution for $\vT_{n,N_n}$ it
is necessary and sufficient that in in (24) the mixing distribution
is the inverse gamma distribution \cite{KorolevSokolov2012}. In
particular, for $\vT_{n,N_n}$ to have the limit multivariate
asymmetric Laplace distribution it is necessary and sufficient that
the limit distribution for $n^{-1}N_n$ is inverse exponential.

\bigskip

\noindent{\bf Acknowledgement.} This research was supported be the
Russian Science Foundation (project 14-11-00364).

\renewcommand{\refname}{References}

\end{document}